# CLOSURES OF EXPONENTIAL FAMILIES[1]

By Imre Csiszár and František Matúš

*Hungarian Academy of Sciences and Academy of Sciences of the Czech Republic*

The variation distance closure of an exponential family with a convex set of canonical parameters is described, assuming no regularity conditions. The tools are the concepts of convex core of a measure and extension of an exponential family, introduced previously by the authors, and a new concept of accessible faces of a convex set. Two other closures related to the information divergence are also characterized.

**1. Introduction.** Exponential families of probability measures (p.m.'s) include many of the parametric families frequently used in statistics, probability and information theory. Their mathematical theory has been worked out to a considerable extent [1, 2, 3, 11]. Although limiting considerations are important and do appear in the literature, less attention has been paid to determining closures of exponential families.

For families supported by a finite or countable set, closures were considered in [1], pages 154–156, and [2], pages 191–201, respectively, the latter with regularity conditions. In the general case, different closure concepts come into account. Our main result, Theorem 2 in Section 3, determines the closure in variation distance (variation closure) of a full exponential family and, more generally, of any subfamily with a convex set of canonical parameters. Weak closures appear much harder to describe in general, but Theorem 1 in Section 3 is a step in that direction.

Received October 2003; revised December 2003.

[1]Supported by the Hungarian National Foundation for Scientific Research under Grants T 26041, T 32323 and TS 40719, and by the Grant Agency of the Academy of Sciences of the Czech Republic under Grant A1075104.

Other closure concepts are based on Kullback–Leibler $I$-divergence (information divergence or relative entropy)

$$D(P\|Q) \triangleq \begin{cases} \int \ln \frac{dP}{dQ}\, dP, & \text{if } P \ll Q, \\ +\infty, & \text{otherwise.} \end{cases}$$

With the terminology of [7], these are the $I$-closure and reverse $I$-closure ($rI$-closure); early work related to the latter appeared in [3]. The $I$-closure of a convex set $\mathcal{S}$ of p.m.'s is relevant, for example, in large deviations theory, where the conditional limit theorem for i.i.d. sequences on the condition that the empirical distribution belongs to $\mathcal{S}$ involves the "generalized $I$-projection" to $\mathcal{S}$ which is in the $I$-closure of $\mathcal{S}$; see [4]. In the context of exponential families, the $rI$-closure is of major statistical interest; for example, when the likelihood function is bounded but its maximum is not attained, a "generalized maximum likelihood estimate" can be introduced as a p.m. that belongs to the $rI$-closure; see [7].

Formally, the variation closure $\mathrm{cl}_v(\mathcal{S})$, respectively, the $I$-closure $\mathrm{cl}_I(\mathcal{S})$ and the $rI$-closure $\mathrm{cl}_{rI}(\mathcal{S})$ of a set $\mathcal{S}$ of p.m.'s on a given measurable space, consists of all p.m.'s $P$ to which there exists a sequence $Q_n$ in $\mathcal{S}$ such that the total variation $|P - Q_n|$, respectively, the $I$-divergence $D(Q_n\|P)$ or $D(P\|Q_n)$, goes to zero as $n \to \infty$. The Pinsker inequality $|P - Q|^2 \leq 2D(P\|Q)$ implies that both $\mathrm{cl}_I(\mathcal{S})$ and $\mathrm{cl}_{rI}(\mathcal{S})$ are contained in $\mathrm{cl}_v(\mathcal{S})$. For exponential families, the last inclusion gives a good approximation to $\mathrm{cl}_{rI}(\mathcal{S})$, for example, all p.m.'s in $\mathrm{cl}_v(\mathcal{S})$ with mean belong to $\mathrm{cl}_{rI}(\mathcal{S})$. This is one motivation for our study of variation closures, in addition to intrinsic mathematical interest. Theorem 3 in Section 3 characterizes those p.m.'s in the variation closure that belong also to the $rI$-closure. The $I$-closure is much easier to describe than the other closures (see Corollary 2), in particular, full exponential families are $I$-closed. It should be mentioned that the $I$- and $rI$-closures are not topological closure operations because they are not idempotent. An example of an exponential family $\mathcal{E}$ with $\mathrm{cl}_{rI}(\mathrm{cl}_{rI}(\mathcal{E}))$ strictly larger than $\mathrm{cl}_{rI}(\mathcal{E})$ is given in [8]. On the other hand, the $I$- and $rI$-closures are sequential closures in suitable topologies; see [9].

Our attention is focused on exponential families that consist of p.m.'s on $\mathbb{R}^d$ and have a canonical statistic equal to the identity mapping. Clearly, determining their variation, $I$- or $rI$-closures, the same problems are solved for general exponential families of p.m.'s on any measurable space, with $d$-dimensional canonical statistic, by mapping the family to one on $\mathbb{R}^d$ via the canonical statistic.

A crucial construction is that of the extension $\mathrm{ext}(\mathcal{E})$ of a full exponential family $\mathcal{E}$, introduced by the authors [5, 7] based on their concept of the convex core of a measure on $\mathbb{R}^d$ [6]; see the definitions in Section 2. The inclusion $\mathcal{E} \subseteq \mathrm{ext}(\mathcal{E})$ is strict unless no nontrivial supporting hyperplane of the (common) convex support of the p.m.'s in $\mathcal{E}$ has positive probability under these p.m.'s, by [7], Remark 3. By Lemma 6 below, the variation closure of $\mathcal{E}$ is contained in $\mathrm{ext}(\mathcal{E})$. A stronger result announced in [5], the variation closedness of $\mathrm{ext}(\mathcal{E})$, follows as Corollary 3 from Theorem 4 that deals with variation convergent sequences in $\mathrm{ext}(\mathcal{E})$.

The inclusion $\mathrm{cl}_v(\mathcal{E}) \subseteq \mathrm{ext}(\mathcal{E})$ implies that if the subset $\mathrm{cl}_{rI}(\mathcal{E})$ of $\mathrm{cl}_v(\mathcal{E})$ is equal to $\mathrm{ext}(\mathcal{E})$ (e.g., if the domain of canonical parameters is the whole $\mathbb{R}^d$; see [7], Lemma 6(ii)), then also $\mathrm{cl}_v(\mathcal{E}) = \mathrm{ext}(\mathcal{E})$. Moreover, since $\mathcal{E}$ is $rI$-closed if and only if $\mathcal{E} = \mathrm{ext}(\mathcal{E})$ ([7], Corollary 2), the last condition is necessary and sufficient also for the variation-closedness of $\mathcal{E}$.

The cases just mentioned, although frequent in practice, are of secondary interest for our purposes. This paper is primarily devoted to the general case when all common regularity conditions are absent, although the assumption of steepness or even regularity (see [1], pages 116 and 117) would not lead to significant simplifications. The typical situation we have in mind is when the p.m.'s in $\mathcal{E}$ have both discrete and continuous components.

## 2. Preliminaries.

2.1. *Convex sets and faces.* The closure and affine hull of a set $B \subseteq \mathbb{R}^d$ are denoted $\mathrm{cl}(B)$ and $\mathrm{aff}(B)$, and the relative interior [interior in the relative topology of $\mathrm{aff}(B)$] is denoted $\mathrm{ri}(B)$. The linear subspace of $\mathbb{R}^d$ obtained by shifting $\mathrm{aff}(B)$ to contain the origin is denoted $\mathrm{lin}(B)$. Orthogonal projections to subspaces of the form $\mathrm{lin}(C)$, where $C \subseteq \mathbb{R}^d$ is a convex set, are often needed in the sequel; they are denoted briefly as $\pi_C$ rather than $\pi_{\mathrm{lin}(C)}$.

A *face* of a nonempty convex set $C \subseteq \mathbb{R}^d$ is a nonempty convex subset $F$ of $C$ such that whenever $tx + (1-t)y \in F$ for some $x, y$ in $C$ and $0 < t < 1$, then $x, y$ are in $F$. As in [10], but unlike in [12], the empty set is not considered to be a face. The *proper* faces are those different from $C$ and the *exposed* faces are $C$ itself and its intersections with the supporting hyperplanes of $C$. A proper exposed face $F$ of $C$ is thus represented as $F = C \cap \{x : \langle \tau, x - a\rangle = 0\}$, where $a \in C$ and $\tau \in \mathbb{R}^d$ is a unit vector such that $\langle \tau, x - a\rangle \leq 0$ for each $x \in C$. Obviously, there is no loss of generality in assuming $\tau \in \mathrm{lin}(C)$. Such a vector $\tau$ *exposes* $F$ in $C$.

2.2. *Convex support and core.* A *measure* always means a finite Borel measure on $\mathbb{R}^d$. The *convex support* $\mathrm{cs}(\mu)$ and the *convex core* $\mathrm{cc}(\mu)$ of $\mu$ are defined, respectively, as the intersection of those convex closed and convex Borel subsets $C$ of $\mathbb{R}^d$ which have full $\mu$-measure, $\mu(C) = \mu(\mathbb{R}^d)$. While the former is a standard concept, the latter is of recent origin [6]. Let us recall from [6] the key facts

$$\text{(1)} \qquad \mathrm{cs}(\mu) = \mathrm{cl}(\mathrm{cc}(\mu)) \quad \text{and} \quad \mathrm{cc}(\mu^{\mathrm{cl}(F)}) = F \qquad \text{for faces } F \text{ of } \mathrm{cc}(\mu),$$

where the restriction of a measure $\mu$ to a Borel subset $B$ of $\mathbb{R}^d$ is denoted $\mu^B$.

Note that the convex closed set $\mathrm{cs}(\mu)$ is of full $\mu$-measure, but the convex set $\mathrm{cc}(\mu)$, though measurable by [6], Theorem 1, need not be. For brevity, $\mathrm{lin}(\mu)$ is written instead of $\mathrm{lin}(\mathrm{cs}(\mu))$ and similarly with $\mathrm{ri}(\mu)$ and $\mathrm{aff}(\mu)$.

LEMMA 1. *A supporting hyperplane $H$ of $\mathrm{cs}(\mu)$ is of positive $\mu$-measure if and only if $F = H \cap \mathrm{cc}(\mu)$ is nonempty. Moreover, $\mu(H \setminus \mathrm{cl}(F)) = 0$.*

PROOF. Using [6], Lemma 2(ii), $F = \mathrm{cc}(\mu^H)$. This and (1) give $\mathrm{cl}(F) = \mathrm{cs}(\mu^H)$, whence both assertions follow. $\square$

2.3. *Exponential families.* The term *exponential family* without any adjective means below a full exponential family based on a (nonzero) measure $\mu$ on $\mathbb{R}^d$, with a canonical statistic equal to the identity mapping. This family $\mathcal{E} = \mathcal{E}_\mu$ consists of the p.m.'s $Q_\vartheta$ with $\mu$-densities

$$\frac{dQ_\vartheta}{d\mu}(x) \triangleq e^{\langle \vartheta, x\rangle - \Lambda(\vartheta)}, \tag{2}$$

where

$$\Lambda(\vartheta) = \Lambda_\mu(\vartheta) \triangleq \ln \int_{\mathbb{R}^d} e^{\langle \vartheta, x\rangle} \mu(dx)$$

and the canonical parameter $\vartheta$ belongs to $\mathrm{dom}(\Lambda) = \{\vartheta : \Lambda(\vartheta) < \infty\}$. Note that $\mu$ is not uniquely determined by the family $\mathcal{E}$. In particular, any member of $\mathcal{E}$ could play the role of $\mu$; in this paper, however, $\mu$ is regarded as given.

Clearly, if $\vartheta \in \mathrm{dom}(\Lambda)$ and $\theta - \vartheta$ is orthogonal to $\mathrm{lin}(\mu)$ for some $\theta \in \mathbb{R}^d$, then also $\theta \in \mathrm{dom}(\Lambda)$ and $Q_\vartheta = Q_\theta$. A bijective parametrization can be given as

$$\begin{aligned} &\mathcal{E} = \{Q_\vartheta : \vartheta \in \Theta\}, \\ &\text{where } \Theta = \Theta_\mu \triangleq \mathrm{dom}(\Lambda_\mu) \cap \mathrm{lin}(\mu) = \pi_{\mathrm{lin}(\mu)}(\mathrm{dom}(\Lambda_\mu)). \end{aligned} \tag{3}$$

Here, $\Theta$ equals $\mathrm{dom}(\Lambda)$ if and only if $\mathrm{lin}(\mu) = \mathbb{R}^d$ [when (2) is called a minimal representation]. For the purposes below, it is essential not to require that condition and not to require $\mu$ to be a p.m., either.

Of main interest are subfamilies

$$\mathcal{E}_\Xi = \{Q_\vartheta : \vartheta \in \Xi\}, \qquad \Xi \subseteq \Theta,$$

of the full family $\mathcal{E}$, primarily when $\Xi$ is convex. In this case, $\mathcal{E}_\Xi$ is called a *canonically convex exponential family.*

The function $\Lambda$ is known to be convex and lower semicontinuous, thus continuous on closed segments contained in $\mathrm{dom}(\Lambda)$. The following lemma is a minor improvement of Lemma 23.5 in [3].

LEMMA 2. *Let $\vartheta_n$ be a sequence in* $\operatorname{dom}(\Lambda)$ *that converges to some* $\vartheta \in \mathbb{R}^d$.

(i) *If $Q_{\vartheta_n}$ converges weakly, then $\vartheta \in \operatorname{dom}(\Lambda)$ and $\Lambda(\vartheta_n)$ converges to $\Lambda(\vartheta)$.*

(ii) *If $\vartheta \in \operatorname{dom}(\Lambda)$ and $\Lambda(\vartheta_n) \to \Lambda(\vartheta)$, then $Q_{\vartheta_n}$ converges to $Q_\vartheta$ in variation.*

PROOF. (i) If $Q_{\vartheta_n}$ converges weakly to some p.m. $P$, then for each continuity set $B$ of $P$,

$$P(B) = \lim_{n\to\infty} Q_{\vartheta_n}(B) = \lim_{n\to\infty} \exp(-\Lambda(\vartheta_n)) \int_B \exp(\langle \vartheta_n, x\rangle)\mu(dx),$$

where, if $B$ is compact,

$$\int_B \exp(\langle \vartheta_n, x\rangle)\mu(dx) \to \int_B \exp(\langle \vartheta, x\rangle)\mu(dx)$$

by dominated convergence. If also $P(B) > 0$, then $\exp(-\Lambda(\vartheta_n))$ converges to a positive number $c$. Hence

$$P(B) = c \cdot \int_B e^{\langle \vartheta, x\rangle}\mu(dx)$$

for each compact continuity set $B$ of $P$ and, consequently, for all Borel sets $B$. When $B = \mathbb{R}^d$, it follows that $c = e^{-\Lambda(\vartheta)}$. Hence, $\vartheta \in \operatorname{dom}(\Lambda)$ and $\Lambda(\vartheta_n)$ converges to $\Lambda(\vartheta)$.

(ii) Under the assumptions, the $\mu$-densities $\exp(\langle \vartheta_n, x\rangle - \Lambda(\vartheta_n))$ of $Q_{\vartheta_n}$ converge to the $\mu$ density $\exp(\langle \vartheta, x\rangle - \Lambda(\vartheta))$ of $Q_\vartheta$ pointwise, which is known to imply $Q_{\vartheta_n} \to Q_\vartheta$ in variation. □

2.4. *Extensions of exponential families.* The restriction of $\mu$ to the closure of a face $F$ of $\operatorname{cc}(\mu)$ is a nonzero measure by (1). The exponential family based on this restriction $\mu^{\operatorname{cl}(F)}$ is denoted $\mathcal{E}^F$. It consists of the p.m.'s $Q_{F,\vartheta}$ defined as in (2) with $\mu$ and $\Lambda$ replaced with $\mu^{\operatorname{cl}(F)}$ and

$$\Lambda_F(\vartheta) \triangleq \ln \int_{\operatorname{cl}(F)} e^{\langle \vartheta, x\rangle}\mu(dx).$$

Obviously, $\operatorname{dom}(\Lambda) \subseteq \operatorname{dom}(\Lambda_F)$. The family $\mathcal{E}^F$ is bijectively parametrized as

$$\mathcal{E}^F = \{Q_{F,\vartheta} : \vartheta \in \Theta_F\}, \qquad \text{where } \Theta_F \triangleq \operatorname{dom}(\Lambda_F) \cap \operatorname{lin}(F) = \pi_F(\operatorname{dom}(\Lambda_F)),$$

similarly to (3), since $\operatorname{lin}(\mu^{\operatorname{cl}(F)}) = \operatorname{lin}(F)$ by (1). For each $Q_\vartheta \in \mathcal{E}$ with $\vartheta \in \operatorname{dom}(\Lambda)$, its conditioning $Q_\vartheta(\cdot|\operatorname{cl}(F))$, equal to the restriction $Q_\vartheta^{\operatorname{cl}(F)}$ divided by $Q_\vartheta(\operatorname{cl}(F))$, belongs to $\mathcal{E}^F$. The simple fact that $Q_\vartheta(\cdot|\operatorname{cl}(F))$ coincides with the p.m. $Q_{F,\theta}$, where $\theta = \pi_F(\vartheta)$ is in $\Theta_F$, is repeatedly used in

the sequel. These conditionings of the p.m.'s in $\mathcal{E}$ exhaust $\mathcal{E}^F$ if and only if the inclusion $\pi_F(\Theta) \subseteq \Theta_F$ holds with equality.

The *extension* $\mathrm{ext}(\mathcal{E})$ of an exponential family $\mathcal{E} = \mathcal{E}_\mu$ is the union of the families $\mathcal{E}^F$ over all faces $F$ of $\mathrm{cc}(\mu)$. Each $\mathcal{E}^F$ is called a *component* of $\mathrm{ext}(\mathcal{E})$. A similar construction of a "boundary at infinity" of $\mathcal{E}$ which uses faces of $\mathrm{cs}(\mu)$ rather than of $\mathrm{cc}(\mu)$ was proposed earlier [3]. Some crucial assertions on an exponential family completed by its "boundary" were found to be erroneous, but their analogues for $\mathrm{ext}(\mathcal{E})$ were found valid in [5, 7]. The reason is that $\mathrm{ext}(\mathcal{E})$ may be strictly larger than $\mathcal{E}$ completed by its "boundary": By [6], Lemma 11, the latter consists of those components $\mathcal{E}^F$ that correspond to the proper exposed faces $F$ of $\mathrm{cc}(\mu)$.

2.5. *Accessible faces.* For any face $F$ of a convex set $C \subseteq \mathbb{R}^d$ there exists a chain

$$C = F_0 \supset F_1 \supset \cdots \supset F_m = F,$$

not necessarily unique, such that $F_i$ is an exposed face of $F_{i-1}$, $1 \le i \le m$. If for every $1 \le i \le m$ a unit vector $\tau_i \in \mathrm{lin}(F_{i-1})$ exposes $F_i$ in $F_{i-1}$, then $\tau_1, \ldots, \tau_m$ is called an *access sequence* to the face $F$ of $C$; the access sequence to $F = C$ is empty. Since $\tau_i \in \mathrm{lin}(F_{i-1})$ is orthogonal to $\mathrm{lin}(F_i)$, the vectors of any nonempty access sequence are orthonormal.

Let $C$ and $\Xi$ be two nonempty convex subsets of $\mathbb{R}^d$. For our main result, where the role of $C$ is played by $\mathrm{cc}(\mu)$ and the role of $\Xi$ is played by a convex subset of $\Theta_\mu$, a new concept of $\Xi$-accessible faces of $C$ is suitable. This concept involves a constraint on access sequences in terms of recession cones of projections of $\mathrm{ri}(\Xi)$. Recall that the *recession cone* of a convex set $\Xi \subseteq \mathbb{R}^d$ is

$$\mathrm{rec}(\Xi) = \{\tau : \vartheta + t\tau \in \Xi \text{ for all } \vartheta \in \Xi, t \ge 0\}.$$

By [12], Theorem 8.2 and Corollary 8.3.1, $\mathrm{rec}(\mathrm{ri}(\Xi)) = \mathrm{rec}(\mathrm{cl}(\Xi))$, and this is a closed cone that contains $\mathrm{rec}(\Xi)$. Now, a face $F$ of the convex set $C$ is $\Xi$-*accessible* if an access sequence $\tau_1, \ldots, \tau_m$ to $F$ exists such that

$$\tau_i \in \mathrm{rec}(\pi_{F_{i-1}}(\mathrm{ri}(\Xi))), \qquad 1 \le i \le m. \tag{4}$$

An access sequence to $F$ that satisfies (4) is called *adapted* to $\Xi$. It may seem artificial that these notions depend on $\Xi$ only through its relative interior, but if ri were omitted in (4), some later assertions would not hold; see Example 3 in Section 3. Note that the empty sequence is trivially adapted; thus $C$ itself is always a $\Xi$-accessible face of $C$.

Lemma 3. *If $\Xi \subseteq \mathrm{lin}(C)$, an access sequence $\tau_1, \ldots, \tau_m$ to a proper face $F$ of $C$ is adapted to $\Xi$ if and only if $\tau_1 \in \mathrm{rec}(\mathrm{ri}(\Xi))$ and for the face $F_1$ of $C$ exposed by $\tau_1$ the access sequence $\tau_2, \ldots, \tau_m$ to the face $F$ of $F_1$ is adapted to $\pi_{F_1}(\Xi)$.*

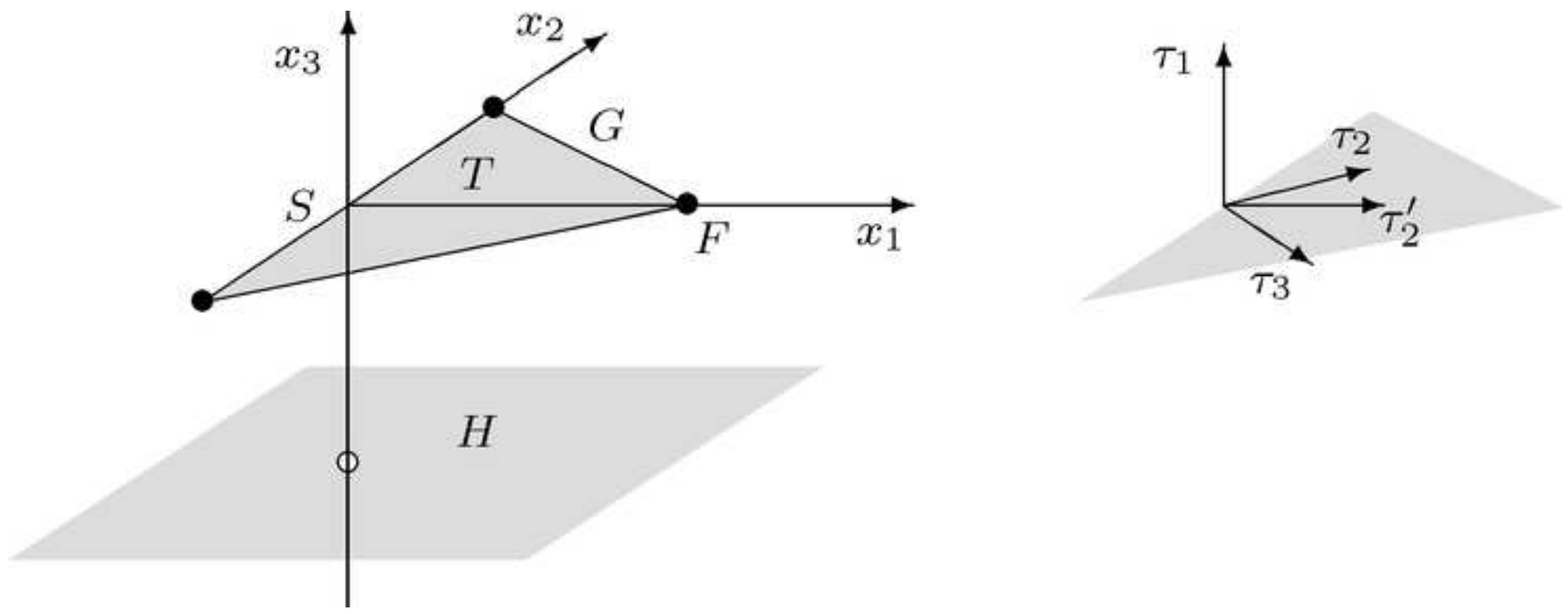

FIG. 1.

PROOF. By the hypotheses $\Xi \subseteq \text{lin}(C)$, the set $\pi_{F_0}(\text{ri}(\Xi))$ in the condition (4) for $i = 1$ is equal to $\text{ri}(\Xi)$. In the conditions for $2 \le i \le m$, the sets $\pi_{F_{i-1}}(\text{ri}(\Xi))$ are equal to the sets $\pi_{F_{i-1}}(\text{ri}(\pi_{F_1}(\Xi)))$ that appear in the analogue of (4) for the adaptedness of $\tau_2, \dots, \tau_m$ to $\pi_{F_1}(\Xi)$, since the operation ri interchanges with orthogonal projections ([12], Theorem 6.6), and $\pi_{F_{i-1}} \pi_{F_1} = \pi_{F_{i-1}}$ if $i \ge 2$. $\square$

EXAMPLE 1 (Figure 1). Let $C \subset \mathbb{R}^3$ be the convex hull of the union of the plane $H = \mathbb{R}^2 \times \{-1\}$ and the triangle $T$ with vertices $(1,0,0)$ and $(0,\pm 1,0)$, and let $\Xi$ consist of those $\vartheta \in \mathbb{R}^3$ whose second coordinate is strictly between $-1$ and 1. In this example, $\Xi$-accessibility of the nine faces of $C$ is discussed. The proper exposed faces of $C$ are $H$ and $T$. Of the six nonexposed faces of $C$, equal to faces of $T$, consider

$$F = \{(1,0,0)\}, \qquad G = \{(t, 1-t, 0) : 0 \le t \le 1\}, \qquad S = \{(0,t,0) : |t| \le 1\}.$$

Since $\Xi$ is open, the relative interiors in (4) can be ignored. Note that the recession cone of $\pi_C(\Xi) = \Xi$ is $\mathbb{R} \times \{0\} \times \mathbb{R}$ and the recession cone of $\pi_T(\Xi)$ is $\mathbb{R} \times \{(0,0)\}$. Let

$$\tau_1 = (0,0,1), \qquad \tau_2 = (\tfrac{1}{\sqrt{2}}, \tfrac{1}{\sqrt{2}}, 0), \qquad \tau_3 = (\tfrac{1}{\sqrt{2}}, -\tfrac{1}{\sqrt{2}}, 0)$$

and

$$\tau_2' = (1,0,0).$$

Since $\text{rec}(\Xi)$ contains both $-\tau_1$ and $\tau_1$, the faces $H$ and $T$ are $\Xi$-accessible. Both $\tau_1, \tau_2, \tau_3$ and $\tau_1, \tau_2'$ are access sequences to $F$, with corresponding chains $C \supset T \supset G \supset F$ and $C \supset T \supset F$, respectively. Since $\text{rec}(\pi_T(\Xi))$ contains $\tau_2'$ but not $\tau_2$, the access sequence $\tau_1, \tau_2'$ is adapted to $\Xi$, whereas $\tau_1, \tau_2, \tau_3$ is not. Due to the former, $F$ is a $\Xi$-accessible face of $C$. On the other hand, $G$ is not $\Xi$-accessible, because the only access sequence $\tau_1, \tau_2$ to $G$, with chain $C \supset T \supset G$ has $\tau_2 \notin \mathrm{rec}(\pi_T(\Xi))$. Similarly, the segment $S$ is $\Xi$-accessible, but its endpoints are not.

2.6. *Partial means.* When studying $rI$-closures of exponential families, p.m.'s that do not have a mean require special attention. The following simple concept is useful: for any p.m. $P$ on $\mathbb{R}^d$, write

$$\mathrm{M}(P) = \{\tau \in \mathbb{R}^d : \langle \tau, \cdot \rangle \text{ is } P\text{-integrable}\}$$

and define the *partial mean* $\mathrm{m}(P)$ as the unique element of the linear space $\mathrm{M}(P)$ with

$$\int_{\mathbb{R}^d} \langle \tau, x \rangle P(dx) = \langle \tau, \mathrm{m}(P) \rangle, \qquad \tau \in \mathrm{M}(P).$$

Note that $\mathrm{M}(P) = \mathbb{R}^d$ if and only if $P$ has a mean, in which case $\mathrm{m}(P)$ equals the mean.

The following lemma is well known, but usually stated in less generality.

Lemma 4. *For $\vartheta \in \mathrm{dom}(\Lambda)$ and a unit vector $\tau$ such that $\vartheta + t\tau \in \mathrm{dom}(\Lambda)$ for some $t > 0$, the integral $\int \langle \tau, x \rangle \, Q_\vartheta(dx)$ exists, either finite or $-\infty$. This integral equals the directional derivative of $\Lambda$ at $\vartheta$ in the direction $\tau$.*

Proof. The directional derivative, that is, the right derivative of the function $t \mapsto \Lambda(\vartheta + t\tau)$ at $t = 0$, equals

$$\frac{1}{\int e^{\langle \vartheta, x \rangle} \mu(dx)} \lim_{t \downarrow 0} \int \frac{e^{t\langle \tau, x \rangle} - 1}{t} e^{\langle \vartheta, x \rangle} \mu(dx).$$

The ratio $(e^{t\langle \tau, x \rangle} - 1)/t$ decreases to $\langle \tau, x \rangle$ when $t$ decreases to zero, and the assertion follows by the monotone convergence. $\square$

The following lemma shows the relevance of partial means for exponential families.

Lemma 5. *The I-divergence $D(Q_\theta \| Q_\vartheta)$ of p.m.'s in $\mathcal{E}$ is finite if and only if $\theta - \vartheta \in \mathrm{M}(Q_\theta)$, in which case*

$$D(Q_\theta \| Q_\vartheta) = \langle \theta - \vartheta, \mathrm{m}(Q_\theta) \rangle - \Lambda(\theta) + \Lambda(\vartheta).$$

Proof. By definition,

$$\begin{aligned} D(Q_\theta \| Q_\vartheta) &= \int_{\mathbb{R}^d} \ln \frac{e^{\langle \theta, x \rangle - \Lambda(\theta)}}{e^{\langle \vartheta, x \rangle - \Lambda(\vartheta)}} Q_\theta(dx) \\ &= \int_{\mathbb{R}^d} \langle \theta - \vartheta, x \rangle Q_\theta(dx) \; - \Lambda(\theta) + \Lambda(\vartheta). \end{aligned}$$

The last integral is finite if and only if $\theta - \vartheta$ belongs to $\mathrm{M}(Q_\theta)$, in which case it equals $\langle \theta - \vartheta, \mathrm{m}(Q_\theta)\rangle$. $\square$

This lemma also is used with $\mathcal{E}$ replaced by a component $\mathcal{E}^F$ of $\mathrm{ext}(\mathcal{E})$, combined with the obvious identity

$$D(P\|Q) = D(P\|Q(\cdot\,|\,\mathrm{cl}(F))) - \ln Q(\mathrm{cl}(F)), \qquad P \in \mathcal{E}^F, Q \in \mathcal{E}. \tag{5}$$

**3. Main results.** Below, $\mathcal{E} = \mathcal{E}_\mu$ always denotes an exponential family in the sense of Section 2.3, $Q_\vartheta$ with $\vartheta \in \Theta = \mathrm{dom}(\Lambda) \cap \mathrm{lin}(\mu)$ denotes a p.m. in $\mathcal{E}$ and $\mathcal{E}_\Xi$ with $\Xi \subseteq \Theta$ denotes a subfamily of $\mathcal{E}$.

THEOREM 1. *If a sequence of p.m.'s $Q_{\vartheta_n}$ with $\vartheta_n$ in $\Theta$ converges weakly to a p.m. $P$, then one of the following two alternatives takes place:*

(i) *The sequence $\vartheta_n$ converges to an element $\vartheta$ of $\Theta$, $P = Q_\vartheta$ and $Q_{\vartheta_n} \to P$ even in variation distance.*

(ii) *The norm of $\vartheta_n$ goes to $\infty$, a proper exposed face $G$ of $\mathrm{cs}(\mu)$ exists such that $P(G) = 1$, and the limit of any convergent subsequence of $\vartheta_n/\|\vartheta_n\|$ exposes such a face of $\mathrm{cs}(\mu)$.*

The proof is given in Section 4. It follows from Theorem 1 that the weak convergence of a sequence $Q_{\vartheta_n}$ to some $Q_\vartheta$ in $\mathcal{E}$ implies $\vartheta_n \to \vartheta$; see [1], Theorem 8.3, for a direct proof. A consequence of this and Lemma 2 is stated for reference purposes.

COROLLARY 1. *For a sequence $\vartheta_n$ in $\Theta$ and $\vartheta \in \Theta$ the following assertions are equivalent.*

(i) *Convergence $Q_{\vartheta_n} \to Q_\vartheta$ weakly.*
(ii) *Convergence $Q_{\vartheta_n} \to Q_\vartheta$ in variation.*
(iii) *Convergences $\vartheta_n \to \vartheta$ and $\Lambda(\vartheta_n) \to \Lambda(\vartheta)$.*

COROLLARY 2. *For any subset $\Xi$ of $\Theta$,*

$$\mathcal{E}_{\mathrm{cl}(\Xi)\cap\Theta} \supseteq \mathrm{cl}_v(\mathcal{E}_\Xi) \cap \mathcal{E} \supseteq \mathrm{cl}_I(\mathcal{E}_\Xi).$$

*When $\Xi$ is convex, the equalities take place.*

PROOF. The first inclusion follows from Corollary 1. For the second one, if $D(Q_{\vartheta_n}\|P) \to 0$, the sequence $Q_{\vartheta_n}$ converges in variation to $P$; thus $\mathrm{cl}_v(\mathcal{E}_\Xi)$ contains $\mathrm{cl}_I(\mathcal{E}_\Xi)$. Moreover, the alternative (i) holds in Theorem 1, since otherwise $D(Q_{\vartheta_n}\|P) = +\infty$ for all $n$. This proves that $\mathcal{E}$ contains $\mathrm{cl}_I(\mathcal{E}_\Xi)$.

Supposing $\Xi$ is convex, let $\vartheta \in \mathrm{cl}(\Xi) \cap \Theta$ and let $\tau$ be a unit vector such that $\vartheta + t\tau$ belongs to $\mathrm{ri}(\Xi)$ for some $t > 0$. It suffices to show that $D(Q_{\vartheta_n}\|Q_\vartheta) \to$ 0 for $\vartheta_n = \vartheta + t_n\tau$ with $t_n$ decreasing to zero. By Lemma 4, $\tau \in \mathrm{M}(Q_{\vartheta_n})$ and then Lemma 5 implies

$$D(Q_{\vartheta_n} \| Q_\vartheta) = t_n \langle \tau, \mathrm{m}(Q_{\vartheta_n})\rangle - \Lambda(\vartheta_n) + \Lambda(\vartheta).$$

Since $\Lambda(\vartheta_n) \to \Lambda(\vartheta)$ and the sequence $\langle \tau, \mathrm{m}(Q_{\vartheta_n})\rangle$ is decreasing by Lemma 4, the claim follows. $\square$

The following technical assertion, which is crucial for Theorem 2, is proved in Section 5. The notion of $\Xi$-accessibility enters the scene.

Lemma 6. *If a sequence of p.m.'s in a canonically convex exponential family $\mathcal{E}_\Xi$ converges in variation distance to a p.m. $P$, then there exists a $\Xi$-accessible face $F$ of $\mathrm{cc}(\mu)$ such that $P$ belongs to $\mathcal{E}^F$.*

Lemma 6 enables us to conclude that the variation closure of a canonically convex family $\mathcal{E}_\Xi$ is contained in the union of those components $\mathcal{E}^F$ of $\mathrm{ext}(\mathcal{E})$ which correspond to the $\Xi$-accessible faces $F$ of $\mathrm{cc}(\mu)$. For our main result, still another tool is needed, namely a special convergence concept: A sequence of p.m.'s $Q_n$ is said to converge *neatly* to a p.m. $P$ if the $P$-dominated component of $Q_n$ has constant $P$-density $c_n$ and $c_n \to 1$ or, equivalently, if $Q_n(\cdot|\,\mathrm{cs}(P))$ equals $P$ when defined and $Q_n(\mathrm{cs}(P)) \to 1$. Obviously, the neat convergence implies variation convergence, and even $rI$ convergence, due to $D(P\|Q_n) = -\ln c_n$.

Lemma 7. *For a convex subset $\Xi$ of $\Theta$ and $\Xi$-accessible face $F$ of $\mathrm{cc}(\mu)$, each p.m. $Q_{F,\theta}$ with $\theta \in \pi_F(\mathrm{ri}(\Xi))$ is the neat limit of a sequence from $\mathcal{E}_{\mathrm{ri}(\Xi)}$.*

This lemma is proved in Section 5.

Theorem 2. *The variation closure of a canonically convex exponential family $\mathcal{E}_\Xi$ is*

$$\mathrm{cl}_v(\mathcal{E}_\Xi) = \bigcup \{Q_{F,\theta} : \theta \in \mathrm{cl}(\pi_F(\Xi)) \cap \Theta_F\},$$

*where the union runs over all $\Xi$-accessible faces $F$ of $\mathrm{cc}(\mu)$.*

Proof. Suppose $P$ is the limit in variation distance of a sequence $Q_n$ in $\mathcal{E}_\Xi$. By Lemma 6, $P \in \mathcal{E}^F$ for a $\Xi$-accessible face $F$ of $\mathrm{cc}(\mu)$. Then $P(\mathrm{cl}(F)) = 1$ implies $Q_n(\mathrm{cl}(F)) \to 1$, and thus the sequence $Q_n(\cdot|\,\mathrm{cl}(F))$ in $\mathcal{E}^F$ also converges to $P$ in variation distance. These conditioned p.m.'s are of form $Q_{F,\theta_n}$ with $\theta_n \in \pi_F(\Xi)$; thus Corollary 2 implies $P = Q_{F,\theta}$ for some $\theta$ in $\mathrm{cl}(\pi_F(\Xi)) \cap \Theta_F$. This proves the inclusion $\subseteq$.

Conversely, if $F$ is a $\Xi$-accessible face of $\mathrm{cc}(\mu)$, Lemma 7 implies that $\mathrm{cl}_v(\mathcal{E}_\Xi)$ contains the p.m.'s $Q_{F,\theta}$ with $\theta \in \pi_F(\mathrm{ri}(\Xi))$. Thus, by Corollary 2 and the convexity of $\pi_F(\mathrm{ri}(\Xi))$,

$$\mathrm{cl}_v(\mathcal{E}_\Xi) \supseteq \{Q_{F,\theta} : \theta \in \mathrm{cl}(\pi_F(\mathrm{ri}(\Xi))) \cap \Theta_F\}.$$

Since ri interchanges with projections, the inclusion $\supseteq$ follows. □

Example 2. Using the notation of Example 1, let $\mu$ be the sum of the measure that gives mass 1 to each vertex of $T$ and of the p.m. on the plane $H$ with density $\exp(-x_1^2)\exp(-|x_2|)$ w.r.t. the Lebesgue measure on $H$. Then $\mathrm{cc}(\mu) = C$ and $\Theta_\mu = \mathrm{dom}(\Lambda_\mu)$ coincides with $\Xi$ of Example 1, defined by restricting the second coordinate to be between $-1$ and 1. By Theorem 2, the variation closure of the full family $\mathcal{E} = \mathcal{E}_\mu$ intersects five out of the nine components of $\mathrm{ext}(\mathcal{E})$. In addition to the full families $\mathcal{E}$, $\mathcal{E}^H$ and $\mathcal{E}^F$, the latter consisting of the point mass at $(1,0,0)$, $\mathrm{cl}_v(\mathcal{E})$ contains also some p.m.'s from $\mathcal{E}^T$ and $\mathcal{E}^S$. Note that $\mathrm{cl}_v(\mathcal{E})$ intersects $\mathcal{E}^F$ but not $\mathcal{E}^G$, although $F \subset G$.

Part (ii) of the following Theorem 3 gives a necessary and sufficient condition for a p.m. $P$ in $\mathrm{cl}_v(\mathcal{E}_\Xi)$ to belong also to $\mathrm{cl}_{rI}(\mathcal{E}_\Xi)$, even when $\Xi \subseteq \Theta$ is not convex. When $\Xi$ is convex, this condition can be effectively verified by Theorem 2. Part (iii) provides a simple sufficient condition which has a direct proof. A trivial consequence of Theorem 3 is that $\mathrm{cl}_{rI}(\mathcal{E}_\Xi)$ contains all p.m.'s in $\mathrm{cl}_v(\mathcal{E}_\Xi)$ that have a mean.

Theorem 3. (i) *Suppose a sequence $Q_{\vartheta_n}$, $\vartheta_n \in \Theta$, of p.m.'s in $\mathcal{E}$ converges in variation to a p.m. $P = Q_{F,\theta}$, $\theta \in \Theta_F$, in a component $\mathcal{E}^F$ of $\mathrm{ext}(\mathcal{E})$. Then $D(P\|Q_{\vartheta_n})$ goes to zero if it is eventually finite, which is equivalent to $\vartheta_n \in \theta + \mathrm{M}(P)$ eventually.*

(ii) *A p.m. $P = Q_{F,\theta}$ belongs to the rI-closure of a subfamily $\mathcal{E}_\Xi$ of $\mathcal{E}$ if and only if it belongs to the variation closure of $\mathcal{E}_{\Xi\cap(\theta+\mathrm{M}(P))}$.*

(iii) *The rI-closure of a canonically convex exponential family $\mathcal{E}_\Xi$ contains all p.m.'s $P$ in the variation closure of $\mathcal{E}_\Xi$ that satisfy $D(P\|Q) < \infty$ for some $Q \in \mathcal{E}_{\mathrm{ri}(\Xi)}$.*

Proof. (i) Consider first the case when $\mathcal{E}^F = \mathcal{E}$; that is, $P = Q_\vartheta$ for some $\vartheta \in \Theta$. Then $\vartheta_n \to \vartheta$ and $\Lambda(\vartheta_n) \to \Lambda(\vartheta)$ by Corollary 1. Thus, the assertion for this case follows by Lemma 5.

When $F$ is a proper face of $\mathrm{cc}(\mu)$, note that the variation convergence assumption implies $Q_{\vartheta_n}(\mathrm{cl}(F)) \to P(\mathrm{cl}(F)) = 1$. Hence, (5) shows that $D(P\|Q_{\vartheta_n})$ goes to zero if and only if $D(P\|Q_{\vartheta_n}(\cdot|\,\mathrm{cl}(F)))$ does. Moreover $Q_{\vartheta_n}(\cdot|\,\mathrm{cl}(F))$, equal to $Q_{F,\theta_n} \in \mathcal{E}^F$ with $\theta_n = \pi_F(\vartheta_n)$, also converges in variation to $P$.

It follows, applying the result in the first case to $\mathcal{E}^F$ in the role of $\mathcal{E}$, that $D(P\|Q_{\vartheta_n}) \to 0$ if it is eventually finite. Also by (5), the finiteness of $D(P\|Q_{\vartheta_n})$ is equivalent to that of $D(P\|Q_{\vartheta_n}(\cdot|\,\mathrm{cl}(F)))$ and, therefore, by Lemma 5, to $\theta-\theta_n \in \mathrm{M}(P)$. Since $\vartheta_n-\theta_n$ is orthogonal to $\mathrm{lin}(F)$, it belongs to $\mathrm{M}(P)$; thus the last condition is equivalent to $\vartheta_n \in \theta+\mathrm{M}(P)$.

(ii) This follows directly from (i).

(iii) Suppose $P = Q_{F,\theta} \in \mathrm{cl}_v(\mathcal{E}_\Xi)$, where $F$ is a $\Xi$-accessible face of $\mathrm{cc}(\mu)$ and $\theta$ belongs to $\mathrm{cl}(\pi_F(\Xi)) \cap \Theta_F$; see Theorem 2. If $D(P\|Q_{\vartheta_0})$ is finite for some $\vartheta_0 \in \mathrm{ri}(\Xi)$, then as in the proof of part (i), $D(P\|Q_{F,\theta_0})$ is also finite, where $\theta_0 = \pi_F(\vartheta_0)$ and $Q_{F,\theta_0} = Q_{\vartheta_0}(\cdot|\,\mathrm{cl}(F))$. Then for $\theta_n = t_n\theta_0 + (1-t_n)\theta$ with $t_n \downarrow 0$, $D(P\|Q_{F,\theta_n})$ is also finite by Lemma 5 and $Q_{F,\theta_n} \to Q_{F,\theta}$ in variation by Corollary 1. It follows by part (i) that $D(P\|Q_{F,\theta_n}) \to 0$.

Since $\theta_0 \in \pi_F(\mathrm{ri}(\Xi)) = \mathrm{ri}(\pi_F(\Xi))$ and $\theta \in \mathrm{cl}(\pi_F(\Xi))$ imply $\theta_n \in \mathrm{ri}(\pi_F(\Xi))$, Lemma 7 gives that each $Q_{F,\theta_n}$ is the neat limit of a sequence in $\mathcal{E}_{\mathrm{ri}(\Xi)}$. Thus to each $\theta_n$ there exists $\vartheta_n \in \mathrm{ri}(\Xi)$ such that $Q_{F,\theta_n} = Q_{\vartheta_n}(\cdot|\,\mathrm{cl}(F))$ and $Q_{\vartheta_n}(\mathrm{cl}(F))$ is arbitrarily close to 1. From this and $D(P\|Q_{F,\theta_n}) \to 0$ the claim $D(P\|Q_{\vartheta_n}) \to 0$ follows by (5). □

The following example illustrates a use of Theorem 3(ii) and Theorem 2 when deciding whether a p.m. belongs to the $rI$-closure of a canonically convex exponential family. It also illustrates why $\mathrm{ri}(\Xi)$ rather than $\Xi$ appears in the definition (4) of $\Xi$-accessibility and in Lemma 7.

Example 3 (Figure 2). Let $\mu$ be the measure on $\mathbb{R}^3$ equal to the sum of the point mass at $(-1,0,0)$, the image $P$ under $t \mapsto (0,t,0)$ of the p.m. with density $\frac{dt}{t^2}$ on the half-line $t>1$ and the image under $t \mapsto (t,t^2,-1)$ of the p.m. with density $\frac{2dt}{t^3}$ on the same half-line. Then

$$\Lambda(\vartheta) = \ln\Bigg[\exp(-\vartheta_1) + \int_1^\infty \exp(\vartheta_2 t)\frac{dt}{t^2} + \exp(-\vartheta_3)\int_1^\infty \exp(\vartheta_1 t + \vartheta_2 t^2)\frac{2\,dt}{t^3}\Bigg], \qquad \vartheta = (\vartheta_1,\vartheta_2,\vartheta_3), \tag{6}$$

$\mathrm{dom}(\Lambda)$ is given by $\vartheta_2<0$ or $\vartheta_2=0$, $\vartheta_1 \le 0$, and $\Theta = \mathrm{dom}(\Lambda)$. Consider $\Xi = \Theta$, thus $\mathcal{E}_\Xi = \mathcal{E}$, and the face $F = \{(0,t,0)\colon t>1\}$ of $\mathrm{cc}(\mu)$. This $F$ is not exposed and the unique access sequence to it is $\tau_1 = (0,0,1)$ and $\tau_2 = (1,0,0)$. This access sequence is adapted to $\Xi$ but it would not be if $\Xi$ rather than $\mathrm{ri}(\Xi)$ had been used in the definition (4). Since $\pi_F(\Xi) = \{(0,t,0)\colon t \le 0\} = \Theta_F$ and $F$ is $\Xi$-accessible, Theorem 2 gives that $\mathcal{E}^F \subseteq \mathrm{cl}_v(\mathcal{E})$. In particular, $\mathrm{cl}_v(\mathcal{E})$ contains $P$ that equals $Q_{F,\theta}$ with $\theta = (0,0,0) \in \pi_F(\Xi)$. On the other hand, as

$$\mathrm{M}(P) = \mathbb{R}\times\{0\}\times\mathbb{R}, \qquad \Xi \cap (\theta + \mathrm{M}(P)) = \{(\vartheta_1,0,\vartheta_3)\colon \vartheta_1 \le 0\}$$

and the unique access sequence to $F$ is not adapted to the latter set, $P$ is not in the variation closure of $\mathcal{E}_{\Xi\cap(\theta+\mathrm{M}(P))}$ by Theorem 2. Consequently $P \notin \mathrm{cl}_{rI}(\mathcal{E})$ by Theorem 3(ii). Thus, $P$ cannot be the neat limit of any sequence in $\mathcal{E}$, showing that Lemma 7 is not valid when $\mathrm{ri}(\Xi)$ is replaced by $\Xi$.

Our final results address variation convergence of arbitrary sequences in $\mathrm{ext}(\mathcal{E})$.

THEOREM 4. *If a sequence $Q_n$ in $\mathrm{ext}(\mathcal{E})$ converges in variation distance to a p.m. $P$, then $P$ belongs to a component $\mathcal{E}^F$ of $\mathrm{ext}(\mathcal{E})$, for sufficiently large $n$ the face $F_n$ of $\mathrm{cc}(\mu)$ with $Q_n \in \mathcal{E}^{F_n}$ contains $F$, and the conditioned p.m.'s $Q_n(\cdot|\,\mathrm{cl}(F))$ belong to $\mathcal{E}^F$ and also converge to $P$ in variation distance.*

Theorem 4 is proved in Section 6.

COROLLARY 3. *The extension $\mathrm{ext}(\mathcal{E})$ of an exponential family $\mathcal{E}$ is variation-closed.*

Corollary 3 strengthens [7], Theorem 2, on the $rI$-closedness of $\mathrm{ext}(\mathcal{E})$. Note that a family $\mathcal{E}$ completed by its "boundary at infinity" in the sense of [3] is not necessarily $rI$-closed, let alone variation-closed, contrary to [3], Lemma 23.7; see [6], Example 3.

COROLLARY 4. *If a sequence $Q_n$ in $\mathrm{ext}(\mathcal{E})$ converges in variation distance to a p.m. $P$ and $D(P\|Q_n)$ is eventually finite, then $D(P\|Q_n) \to 0$.*

The eventual finiteness takes place, in particular, if $P$ has a mean.

PROOF OF COROLLARY 4. By Theorem 4, the variation convergence of $Q_n \in \mathcal{E}^{F_n}$ to $P$ implies $P \in \mathcal{E}^F$, $F \subseteq F_n$ eventually and $Q_n(\cdot|\,\mathrm{cl}(F)) \to P$ in variation. Hence, the proof can be completed similarly to that of Theorem 3(i) using (5) with $\mathcal{E}^{F_n}$ playing the role of $\mathcal{E}$. $\square$

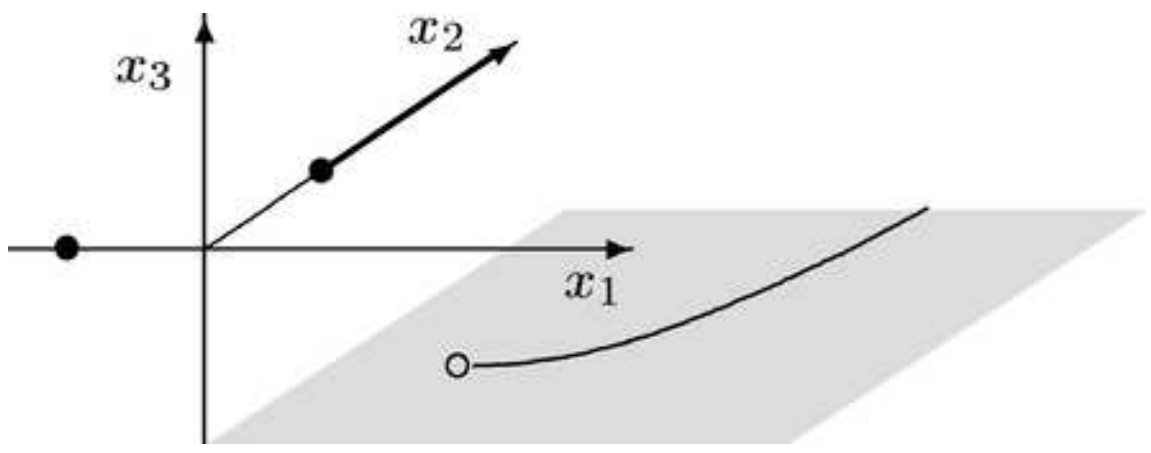

FIG. 2.

**4. Weak convergence in exponential families.** In this section Theorem 1 is proved. For its second alternative, a corollary of the following lemma is needed. Note that

$$\rho(a) \triangleq \inf\{\|y-a\| : y \in \operatorname{aff}(\mu) \setminus \operatorname{cs}(\mu)\}$$

is obviously positive for any element $a$ of $\operatorname{ri}(\mu)$.

LEMMA 8. *For $a \in \operatorname{ri}(\mu)$ and $0 < s < \rho(a)$ there exists a positive constant $C$ such that the inequality*

$$\langle \vartheta, b-a\rangle \geq r\|\vartheta\| - Ce^{-s\|\vartheta\|}[r\|\vartheta\|e^{r\|\vartheta\|}+1] \tag{7}$$

*holds for all $0 < r < s$, $\vartheta \in \Theta$ and $b \in \mathbb{R}^d$ satisfying $\langle \vartheta, b\rangle = \int_{\mathbb{R}^d} \langle \vartheta, x\rangle Q_\vartheta(dx)$.*

In particular, (7) holds when $b$ is the mean of $Q_\vartheta$. What actually is used is the following consequence of Lemma 8.

COROLLARY 5. *Let $\vartheta_n$ be a sequence in $\Theta$ such that each p.m. $Q_{\vartheta_n}$ has a mean $b_n$. If $\|\vartheta_n\| \to \infty$, $\vartheta_n/\|\vartheta_n\| \to \tau$ and $b_n \to b$, then $\langle \tau, b-a\rangle \geq \rho(a)$ for each $a \in \operatorname{ri}(\mu)$.*

PROOF. For $a \in \operatorname{ri}(\mu)$ and $0 < r < s < \rho(a)$, (7) implies

$$\langle \vartheta_n, b_n - a\rangle \geq r\|\vartheta_n\| - C[r\|\vartheta_n\| \exp(-(s-r)\|\vartheta_n\|) + \exp(-s\|\vartheta_n\|)].$$

Dividing both sides by $\|\vartheta_n\|$ and letting $n \to \infty$, it follows that $\langle \tau, b-a\rangle \geq r$ for all $0 < r < s < \rho(a)$. $\square$

For the proof of Lemma 8 the following auxiliary lemma, a simple refinement of known assertions (see [3], Lemma 21.8, and [11], proof of Theorem 3.1), is needed. Denote

$$A^\vartheta_{a,s} \triangleq \{x \in \mathbb{R}^d : \langle \vartheta, x-a\rangle \geq s\|\vartheta\|\}, \qquad \vartheta, a \in \mathbb{R}^d, s > 0,$$

which is a closed half-space of distance $s$ from $a$ when $\vartheta \neq 0$.

LEMMA 9. *If $a \in \operatorname{ri}(\mu)$ and $0 < s < \rho(a)$, then*

$$c_{a,s} \triangleq \inf\{\mu(A^\vartheta_{a,s}) : \vartheta \in \operatorname{lin}(\mu)\}$$

*is positive and*

$$\Lambda(\vartheta) - \langle \vartheta, a\rangle - s\|\vartheta\| \geq \ln c_{a,s}, \qquad \vartheta \in \operatorname{lin}(\mu). \tag{8}$$

Proof. Let $\mu(A_{a,s}^{\vartheta_n}) \downarrow c_{a,s}$ for a sequence $\vartheta_n$ in $\mathrm{lin}(\mu)$ that can be supposed to consist of unit vectors converging to some $\tau \in \mathrm{lin}(\mu)$. Each $x \in \mathbb{R}^d$ with $\langle \tau, x-a\rangle > s$ belongs to $A_{a,s}^{\vartheta_n}$ eventually; thus the open half-space given by the last inequality is covered by the union over all $m$ of the intersections $\bigcap_{n\geq m} A_{a,s}^{\vartheta_n}$. Since the $\mu$-measure of the half-space is positive due to $s < \rho(a)$, one of these intersections has positive $\mu$-measure. Thus, $c_{a,s} > 0$. Finally, the inequalities

$$e^{\Lambda(\vartheta)-\langle\vartheta,a\rangle} = \int_{\mathbb{R}^d} e^{\langle\vartheta,x-a\rangle}\mu(dx) \geq \int_{A_{a,s}^{\vartheta}} e^{\langle\vartheta,x-a\rangle}\mu(dx) \geq e^{s\|\vartheta\|}\cdot\mu(A_{a,s}^{\vartheta})$$

are valid for any $\vartheta \in \mathrm{lin}(\mu)$ and imply (8). □

Proof of Lemma 8. Abbreviate $A_{a,r}^{\vartheta}$ to $A$ and abbreviate its complement to $B$. Then

$$\langle\vartheta, b-a\rangle = \int_{\mathbb{R}^d}\langle\vartheta,x-a\rangle Q_\vartheta(dx) \geq \int_B \langle\vartheta,x-a\rangle Q_\vartheta(dx) + r\|\vartheta\|Q_\vartheta(A)$$

and, in turn,

$$\begin{aligned}\langle\vartheta, b-a\rangle - r\|\vartheta\| &\geq \int_B [\langle\vartheta,x-a\rangle - r\|\vartheta\|]Q_\vartheta(dx)\\ &= e^{\langle\vartheta,a\rangle-\Lambda(\vartheta)}\int_B[\langle\vartheta,x-a\rangle - r\|\vartheta\|]e^{\langle\vartheta,x-a\rangle}\mu(dx).\end{aligned}$$

Since $te^t > -1$ for each $t \in \mathbb{R}$ and $e^{\langle\vartheta,x-a\rangle} < e^{r\|\vartheta\|}$ for $x \in B$, the last integral is bounded below by

$$-[1 + r\|\vartheta\|e^{r\|\vartheta\|}]\mu(B).$$

It follows that

$$\langle\vartheta, b-a\rangle - r\|\vartheta\| \geq -e^{\langle\vartheta,a\rangle-\Lambda(\vartheta)+s\|\vartheta\|}\mu(\mathbb{R}^d)\cdot e^{-s\|\vartheta\|}[r\|\vartheta\|e^{r\|\vartheta\|}+1].$$

This completes the proof on account of (8). [The constant $C$ in (7) can be chosen as $\mu(\mathbb{R}^d)/c_{a,s}$.] □

Proof of Theorem 1. Let a sequence of p.m.'s $Q_{\vartheta_n}$ with $\vartheta_n$ in $\Theta = \mathrm{dom}(\Lambda)\cap\mathrm{lin}(\mu)$ converge weakly to a p.m. $P$. Note that $P(\mathrm{cs}(\mu)) = 1$ because $Q(\mathrm{cs}(\mu)) = 1$ for every $Q \in \mathcal{E}$.

In the case when the sequence $\vartheta_n$ is bounded, let $\vartheta$ be the limit of an arbitrary convergent subsequence. By Lemma 2(i), the weak convergence of $Q_{\vartheta_n}$ along this subsequence implies $\vartheta \in \mathrm{dom}(\Lambda)$, and since clearly $\vartheta \in \mathrm{lin}(\mu)$, also $\vartheta \in \Theta$. Moreover, $P = Q_\vartheta$ by Lemma 2, and since the parametrization (3) of $\mathcal{E}$ is bijective, it follows that each convergent subsequence of $\vartheta_n$ has the same limit. Thus $\vartheta_n$ converges to $\vartheta$. Applying Lemma 2, $\Lambda(\vartheta_n) \to \Lambda(\vartheta)$ follows from the weak convergence of $Q_{\vartheta_n}$ by (i), and $Q_{\vartheta_n} \to P$ in variation by (ii).

In the case when $\|\vartheta_n\| \to \infty$, assume with no loss of generality that $\vartheta_n/\|\vartheta_n\|$ converges to some $\tau$, clearly in $\mathrm{lin}(\mu)$. Suppose first that $\mathrm{cs}(\mu)$ is compact. Then the mean $b_n$ of $Q_{\vartheta_n}$ exists and converges to the mean $b$ of $P$. By Corollary 5, $\langle \tau, b-a\rangle \geq \rho(a)$ for each $a \in \mathrm{ri}(\mu)$ and thus $\langle \tau, x-b\rangle \leq 0$ for each $x \in \mathrm{cs}(\mu)$. Since obviously $b \in \mathrm{cs}(\mu)$, it follows that $H_b = \{x : \langle \tau, x-b\rangle = 0\}$ is a supporting hyperplane of $\mathrm{cs}(\mu)$, and $G = H_b \cap \mathrm{cs}(\mu)$ is a proper face of $\mathrm{cs}(\mu)$, exposed by $\tau$. Then $P(\mathrm{cs}(\mu)) = 1$ implies $P(G) = P(H_b)$ and this equals 1 because the mean $b$ of $P$ belongs to $H_b$.

Turning to the situation when $\mathrm{cs}(\mu)$ is not compact, there exists a continuity set $B$ of $P$ such that $P(B) > 0$ and $\mathrm{lin}(\mu^B) = \mathrm{lin}(\mu)$. Then the conditioned p.m.'s $Q_{\vartheta_n}(\cdot|B)$ belong to the exponential family based on $\mu^B$, the parameters $\vartheta_n$ in $\Theta = \Theta_\mu$ belong also to $\Theta_{\mu^B}$ and $Q_{\vartheta_n}(\cdot|B)$ converges weakly to $P(\cdot|B)$. If, in addition, $B$ is compact, then $P(\cdot|B)$ has a mean $b$ and the result proved above gives that $H_b$ is a nontrivial supporting hyperplane of $\mathrm{cs}(\mu^B)$. Taking another compact continuity set $C$ of $P$ with $C \supseteq B$, let $c$ be the mean of $P(\cdot|C)$ and let $H_c$ be the corresponding supporting hyperplane of $\mathrm{cs}(\mu^C)$. Since $P(H_b|B) = 1$ and $P(H_c|C) = 1$ together with $C \supseteq B$ imply $P(H_c|B) = 1$, the parallel hyperplanes $H_b$ and $H_c$ coincide. This proves that $H_b$ is a nontrivial supporting hyperplane of $\mathrm{cs}(\mu^C)$ satisfying $P(H_b|C) = 1$ for those compact continuity sets $C$ of $P$ which contain $B$. Then $P(H_b) = 1$ and $H_b$ is a nontrivial supporting hyperplane to $\mathrm{cs}(\mu)$, as well, because each $x \in \mathrm{ri}(\mu)$ belongs to $\mathrm{cs}(\mu^C)$ for some $C$ as above. Thus, $G = H_b \cap \mathrm{cs}(\mu)$ is a proper face of $\mathrm{cs}(\mu)$ exposed by $\tau$ and $P(G) = 1$, the same conclusions as when $\mathrm{cs}(\mu)$ was compact.

Finally, $\vartheta_n$ cannot have both a convergent subsequence and another subsequence with norms tending to infinity. Indeed, by the above arguments, the former implies $P \in \mathcal{E}$ and the latter implies $P(G) = 1$ for a proper face of $\mathrm{cs}(\mu)$, a contradiction. □

## 5. Proofs of Lemmas 6 and 7.

Proof of Lemma 6. By induction on the dimension of $\mathrm{cs}(\mu)$. The case of zero dimension is trivial. The induction hypothesis assumes the assertion is true for canonically convex exponential families based on measures whose convex support has smaller dimension than that of the given $\mu$.

Given $Q_{\vartheta_n}$ with $\vartheta_n \in \Xi$ converging in variation distance to $P$, if the alternative (i) in Theorem 1 takes place, then $P \in \mathcal{E}^F$ holds with $F = \mathrm{cc}(\mu)$, obviously a $\Xi$-accessible face of $\mathrm{cc}(\mu)$. Otherwise, by Theorem 1(ii), $\|\vartheta_n\| \to +\infty$ and the limit $\tau_1$ of any convergent subsequence of $\vartheta_n/\|\vartheta_n\|$ exposes a proper face $G$ of $\mathrm{cs}(\mu)$ with $P(G) = 1$. Note that $\vartheta_n \in \Theta \subseteq \mathrm{lin}(\mu)$ implies $\tau_1 \in \mathrm{lin}(\mu)$. The supporting hyperplane $H = \{x : \langle \tau_1, x-a\rangle = 0\}$ that contains $G$ has positive $\mu$-measure, since $P(H) = P(G) = 1$ and $\mu$ dominates the variation limit $P$ of the p.m.'s $Q_{\vartheta_n}$. It follows by Lemma 1 that $F_1 = H \cap \mathrm{cc}(\mu)$ is a face of $\mathrm{cc}(\mu)$, clearly exposed by $\tau_1$, and $\mu(H \setminus \mathrm{cl}(F_1)) = 0$. Thus $P(H) = 1$ implies $P(\mathrm{cl}(F_1)) = 1$. By the variation convergence $Q_{\vartheta_n} \to P$, it follows that $Q_{\vartheta_n}(\mathrm{cl}(F_1)) \to 1$ and the conditioned p.m.'s $Q_{\vartheta_n}(\cdot|\,\mathrm{cl}(F_1))$ also converge to $P$ in variation distance. These p.m.'s belong to $\mathcal{E}^{F_1}$ and thus can be written as $Q_{F_1,\theta_n}$ with $\theta_n = \pi_{F_1}(\vartheta_n) \in \pi_{F_1}(\Xi)$. Hence, by the induction hypothesis applied to the canonically convex exponential family

$$\{Q_{F_1,\theta} \in \mathcal{E}^{F_1} : \theta \in \pi_{F_1}(\Xi)\},$$

their variation limit $P$ belongs to $\mathcal{E}^F$ for a face $F$ of $\mathrm{cc}(\mu^{\mathrm{cl}(F_1)}) = F_1$, which is $\pi_{F_1}(\Xi)$-accessible, that is, an access sequence $\tau_2, \ldots, \tau_m$ to the face $F$ of $F_1$ is adapted to $\pi_{F_1}(\Xi)$. Since $\tau_1$ is the limit of a convergent subsequence of $\vartheta_n/\|\vartheta_n\|$ with $\vartheta_n \in \Xi$, it belongs to $\mathrm{rec}(\mathrm{ri}(\Xi))$ by [12], Theorem 8.2. This and the adaptedness of $\tau_2, \ldots, \tau_m$ to $\pi_{F_1}(\Xi)$ imply by Lemma 3 that the access sequence $\tau_1, \ldots, \tau_m$ to the face $F$ of $\mathrm{cc}(\mu)$ is adapted to $\Xi$. □

The following simple auxiliary assertion resembles [3], Lemma 21.7, and [7], Lemmas 6 and 7(i).

Lemma 10. *If $\tau$ exposes a face $G$ of $\mathrm{cc}(\mu)$, then $\tau$ belongs to the recession cone of $\mathrm{dom}(\Lambda_\mu)$ and for every $\vartheta \in \mathrm{dom}(\Lambda_\mu)$, the sequence $Q_{\vartheta+n\tau}$ converges neatly to $Q_\vartheta(\cdot|\,\mathrm{cl}(G)) \in \mathcal{E}^G$.*

Proof. Since $\tau$ exposes $G$, $\langle \tau, x\rangle \le \langle \tau, a\rangle$ for all $x \in \mathrm{cc}(\mu)$ and $a \in G$, with equality if and only if $x \in G$. Then, for $t \ge 0$ the function $e^{\langle \vartheta + t\tau, x\rangle}$ of $x \in \mathrm{cs}(\mu)$ is bounded above by $e^{\langle \vartheta, x\rangle} \cdot e^{t\langle \tau, a\rangle}$. This implies that $\vartheta + t\tau \in \mathrm{dom}(\Lambda_\mu)$ whenever $\vartheta \in \mathrm{dom}(\Lambda_\mu)$ and $t \ge 0$, proving $\tau \in \mathrm{rec}(\mathrm{dom}(\Lambda_\mu))$.

Knowing from Lemma 1 that $\mu(\mathrm{cl}(G))$ equals the $\mu$-measure of the supporting hyperplane $H = \{x : \langle \tau, x - a\rangle = 0\}$ with $a \in G$, we have for any $\vartheta \in \mathrm{dom}(\Lambda_\mu)$

$$\begin{aligned}
&\Lambda_\mu(\vartheta + n\tau) - \langle \vartheta + n\tau, a\rangle \\
&\qquad = \ln\left[\int_{\mathbb{R}^d} e^{\langle \vartheta + n\tau, x - a\rangle} \mu(dx)\right] \\
&\qquad = \ln\left[\int_{\mathrm{cl}(G)} e^{\langle \vartheta, x - a\rangle} \mu(dx) + \int_{\mathrm{cs}(\mu)\setminus H} e^{\langle \vartheta + n\tau, x - a\rangle} \mu(dx)\right].
\end{aligned}$$

When $n$ tends to infinity, $\Lambda_\mu(\vartheta + n\tau) - \langle \vartheta + n\tau, a\rangle$ decreases to $\Lambda_G(\vartheta) - \langle \vartheta, a\rangle$, since the integral over $\mathrm{cs}(\mu) \setminus H$ decreases to zero by dominated convergence. This fact and $e^{\langle \vartheta + n\tau, x\rangle} = e^{\langle \vartheta, x\rangle + n\langle \tau, a\rangle}$ for $x \in \mathrm{cl}(G)$ imply

$$\begin{aligned}
Q_{\vartheta+n\tau}(\mathrm{cl}(G)) &= e^{\langle \vartheta + n\tau, a\rangle - \Lambda(\vartheta + n\tau)} \int_{\mathrm{cl}(G)} e^{\langle \vartheta, x - a\rangle} \mu(dx) \\
&\to \int_{\mathrm{cl}(G)} e^{\langle \vartheta, x\rangle - \Lambda_G(\vartheta)} \mu(dx) = 1
\end{aligned}$$

and $Q_{\vartheta+n\tau}(\cdot|\,\mathrm{cl}(G)) = Q_\vartheta(\cdot|\,\mathrm{cl}(G))$. Thus the neat convergence follows. □

Proof of Lemma 7. By induction. As in the proof of Lemma 6, the induction hypothesis assumes the assertion is true for exponential families based on measures whose convex support has smaller dimension than that of $\mu$.

The assertion trivially holds if $F = \mathrm{cc}(\mu)$. Thus suppose $F$ is proper face of $\mathrm{cc}(\mu)$ and let $\tau_1, \ldots, \tau_m$ be an access sequence to $F$ adapted to $\Xi$. Let $F_1$ be the face of $\mathrm{cc}(\mu)$ exposed by $\tau_1$. Then, by Lemma 3, $\tau_1 \in \mathrm{rec}(\mathrm{ri}(\Xi))$ and the access sequence $\tau_2, \ldots, \tau_m$ to the face $F$ of $F_1$ is adapted to $\pi_{F_1}(\Xi)$. Thus $F$ is a $\pi_{F_1}(\Xi)$-accessible face of $F_1$.

To prove that for $\vartheta \in \mathrm{ri}(\Xi)$ there exists a sequence of p.m.'s in $\mathcal{E}_{\mathrm{ri}(\Xi)}$ that converges neatly to $Q_\vartheta(\cdot|\,\mathrm{cl}(F))$, apply the induction hypothesis to the exponential family $\mathcal{E}^{F_1}$ based on $\mu^{\mathrm{cl}(F_1)}$ with convex core $F_1$, to the $\pi_{F_1}(\Xi)$-accessible face $F$ of $F_1$ and to $\theta = \pi_{F_1}(\vartheta)$ in $\pi_{F_1}(\mathrm{ri}(\Xi)) = \mathrm{ri}(\pi_{F_1}(\Xi))$. It follows that some sequence $Q_{F_1,\theta_n}$ in $\mathcal{E}^{F_1}$ with $\theta_n \in \pi_{F_1}(\mathrm{ri}(\Xi))$ converges neatly to the conditioning on $\mathrm{cl}(F)$ of $Q_{F_1,\theta}$, which equals $Q_\vartheta(\cdot|\,\mathrm{cl}(F))$ since $\theta = \pi_{F_1}(\vartheta)$. Here, on account of $\theta_n \in \pi_{F_1}(\mathrm{ri}(\Xi))$, the p.m. $Q_{F_1,\theta_n}$ equals $Q_{\vartheta_n}(\cdot|\,\mathrm{cl}(F_1))$ for some $\vartheta_n \in \mathrm{ri}(\Xi)$.

Since $\tau_1 \in \mathrm{rec}(\mathrm{ri}(\Xi))$, Lemma 10 gives that the sequence $Q_{\vartheta_n+k\tau_1}$ converges neatly to $Q_{\vartheta_n}(\cdot|\,\mathrm{cl}(F_1)) = Q_{F_1,\theta_n}$, where each $Q_{\vartheta_n+k\tau_1}$ is in $\mathcal{E}_{\mathrm{ri}(\Xi)}$ due to $\tau_1 \in \mathrm{rec}(\mathrm{ri}(\Xi))$. The last assertion and the neat convergence of $Q_{F_1,\theta_n}$ to $Q_\vartheta(\cdot|\,\mathrm{cl}(F))$ imply that for a suitable sequence $k_n \to \infty$, the p.m.'s $Q_{\vartheta_n+k_n\tau_1}$ in $\mathcal{E}_{\mathrm{ri}(\Xi)}$ converge neatly to $Q_\vartheta(\cdot|\,\mathrm{cl}(F))$. □

**6. Variation convergence in ext($\mathcal{E}$).** In this section Theorem 4 is proved. An auxiliary lemma is sent forward.

Lemma 11. *If $\mu$ dominates a p.m. $P$, then:*

(i) *There exists a face $F$ of $\mathrm{cc}(\mu)$ with $P(\mathrm{cl}(F)) = 1$ such that all faces $G$ of $\mathrm{cc}(\mu)$ with $P(\mathrm{cl}(G)) = 1$ contain $F$.*

(ii) *If $P(\mathrm{cl}(F_n)) \to 1$ for a sequence $F_n$ of proper faces of $\mathrm{cc}(\mu)$, then the face $F$ of* (i) *is proper.*

Proof. (i) The closure $\mathrm{cs}(\mu)$ of $\mathrm{cc}(\mu)$ has full $\mu$-measure, hence also full $P$-measure due to domination. Thus, the face $G = \mathrm{cc}(\mu)$ of $\mathrm{cc}(\mu)$ satisfies $P(\mathrm{cl}(G)) = 1$. Consider any face $G$ with the last property and let $F$ be a face with that property whose dimension is smallest. Then both $\mu^{\mathrm{cl}(G)}$ and $\mu^{\mathrm{cl}(F)}$ dominate $P$, hence so does also the restriction of $\mu$ to $\mathrm{cl}(F) \cap \mathrm{cl}(G)$. By [6], Corollary 4, this intersection has the same $\mu$-measure as $\mathrm{cl}(F \cap G)$. Therefore, the restriction of $\mu$ to $\mathrm{cl}(F \cap G)$ dominates $P$ and thus $P(\mathrm{cl}(F \cap G)) = 1$. The minimality of the dimension of $F$ implies $F \subseteq G$.

(ii) The proper faces $F_n$ in the hypotheses can be supposed to be exposed. Thus let a unit vector $\tau_n$ from $\mathrm{lin}(\mu)$ expose $F_n$ of $\mathrm{cc}(\mu)$. Then for $a \in \mathrm{ri}(\mu)$ the closed half-space $\{x : \langle \tau_n, x-a\rangle \le 0\}$ is disjoint with $\mathrm{cl}(F_n)$; thus its $P$-measure is at most $1 - P(\mathrm{cl}(F_n))$. It can be assumed that $\tau_n \to \tau$ and then, as in the proof of Lemma 9,

$$\{x : \langle \tau, x-a\rangle < 0\} \subseteq \bigcup_{m\ge 1} \bigcap_{n\ge m} \{x : \langle \tau_n, x-a\rangle \le 0\}.$$

Since $P(\{x : \langle \tau_n, x-a\rangle \le 0\}) \le 1 - P(\mathrm{cl}(F_n))$ and $P(\mathrm{cl}(F_n)) \to 1$, the open half-space on the left-hand side has $P$-measure zero whenever $a \in \mathrm{ri}(\mu)$. Hence, on account of $P(\mathrm{cs}(\mu)) = 1$, $\tau$ exposes a proper face of $\mathrm{cs}(\mu)$ that has full $P$-measure. Thus, there exists a nontrivial supporting hyperplane $H$ of $\mathrm{cs}(\mu)$ with $\mu(H) > 0$. By Lemma 1, $G = H \cap \mathrm{cc}(\mu)$ is a proper face of $\mathrm{cc}(\mu)$ and $\mu(H \setminus \mathrm{cl}(G)) = 0$. It follows that $P(\mathrm{cl}(G)) = P(H) = 1$. Hence, $G$ contains the face $F$ of (i) which implies that $F$ is proper. □

Proof of Theorem 4. The variation limit $P$ of p.m.'s $Q_n$ in $\mathrm{ext}(\mathcal{E})$ is obviously dominated by $\mu$; thus Lemma 11(i) applies to this $P$. Let $F$ be the smallest face of $\mathrm{cc}(\mu)$ with closure of full $P$-measure. The variation convergence $Q_n \to P$ implies $Q_n(\mathrm{cl}(F)) \to P(\mathrm{cl}(F)) = 1$, and $Q_n \in \mathcal{E}^{F_n}$ implies $Q_n(\mathrm{cl}(F_n)) \to 1$. Since $Q_n(\mathrm{cl}(F_n) \cap \mathrm{cl}(F))$ equals $Q_n(\mathrm{cl}(F_n \cap F))$ by [6], Corollary 4, it follows that $Q_n(\mathrm{cl}(F_n \cap F)) \to 1$. Thus, again by the variation convergence, also $P(\mathrm{cl}(F_n \cap F)) \to 1$. If a subsequence of $F \cap F_n$ consisted of proper faces of $F$, the last limit relationship would imply, by Lemma 11(ii), applied to $\mu^{\mathrm{cl}(F)}$ in the role of $\mu$, the existence of a proper face of $F = \mathrm{cc}(\mu^{\mathrm{cl}(F)})$ with closure of full $P$-measure, a contradiction to the choice of $F$. This proves that $F_n$ eventually contains $F$. The last inclusion implies that the conditioning $Q_n(\cdot|\,\mathrm{cl}(F))$ of $Q_n \in \mathcal{E}^{F_n}$ belongs to $\mathcal{E}^F$, and since $Q_n(\mathrm{cl}(F)) \to 1$, these conditionings also converge to $P$ in variation distance. Finally, applying Theorem 1 to the p.m.'s $Q_n(\cdot|\,\mathrm{cl}(F))$ in $\mathcal{E}^F$, the alternative (ii) is ruled out by $F$ not having proper faces with closure of full $P$-measure, and it follows that $P \in \mathcal{E}^F$. □

A. Rényi Institute of Mathematics
Hungarian Academy of Sciences
H-1364 Budapest
P.O. Box 127
Hungary
e-mail: csiszar@renyi.hu

Institute of Information Theory and Automation
Academy of Sciences of the Czech Republic
Pod vodárenskou věží 4
182 08 Prague
Czech Republic
e-mail: matus@utia.cas.cz